\def\timestamp{%
Time-stamp: <simon.tex: Thursday 03-09-2020 at 17:16:13 (cest)>}
\def\stripname Time-stamp: <#1 #2>{#2}
\edef\filedate{\expandafter\stripname\timestamp}
\documentclass[a4paper]{amsart}

\usepackage{amssymb}
\usepackage{graphicx}

\let\germ\mathfrak
\let\axiom\mathsf
\newcommand\CH{\axiom{CH}}
\newcommand\ZFC{\axiom{ZFC}}
\newcommand\PFA{\axiom{PFA}}
\newcommand\bee{\mathfrak{b}}
\newcommand\cee{\mathfrak{c}}
\newcommand\dee{\mathfrak{d}}
\newcommand\haa{\mathfrak{h}}
\newcommand\enn{\mathfrak{n}}
\newcommand\pee{\mathfrak{p}}
\newcommand\tee{\mathfrak{t}}
\DeclareMathSymbol\N 0{AMSb}{`N}
\DeclareMathSymbol\Q 0{AMSb}{`Q}
\DeclareMathSymbol\Z 0{AMSb}{`Z}
\newcommand\betaN{\beta\N}
\newcommand\betaNminN{\betaN\setminus\N}
\newcommand\Nstar{\N^*}
\newcommand\fin{\mathit{fin}}
\newcommand\calA{\mathcal{A}}
\newcommand\calI{\mathcal{I}}

\newcommand\pow{\mathcal{P}}
\newcommand\Pkappamodsmall{\pow(\kappa)/[\kappa]^{<\kappa}}
\usepackage{amsrefs}
\theoremstyle{definition}
\newtheorem{question}{Question}[section]

\begin{document}

\title{Petr Simon (1944-2018)}

\author{K. P. Hart}

\address{TU Delft, Delft, The Netherlands}
\email{k.p.hart@tudelft.nl}

\author{M. Hru\v{s}\'ak}
\address{Centro de Ciencias Matem\'aticas, UNAM, Morelia, Mexico}
\email{michael@matmor.unam.mx}

\author{J. L. Verner}
\address{Faculty of Philosophy, Charles University, Prague, Czech Republic}
\email{jonathan.verner@ff.cuni.cz}

\begin{abstract} 
This article is a reflection on the mathematical legacy of Professor 
Petr Simon. 
\end{abstract}

\date{\filedate}



\maketitle

\begin{figure}[ht]
\begin{center} 
\includegraphics[width=.7\hsize]{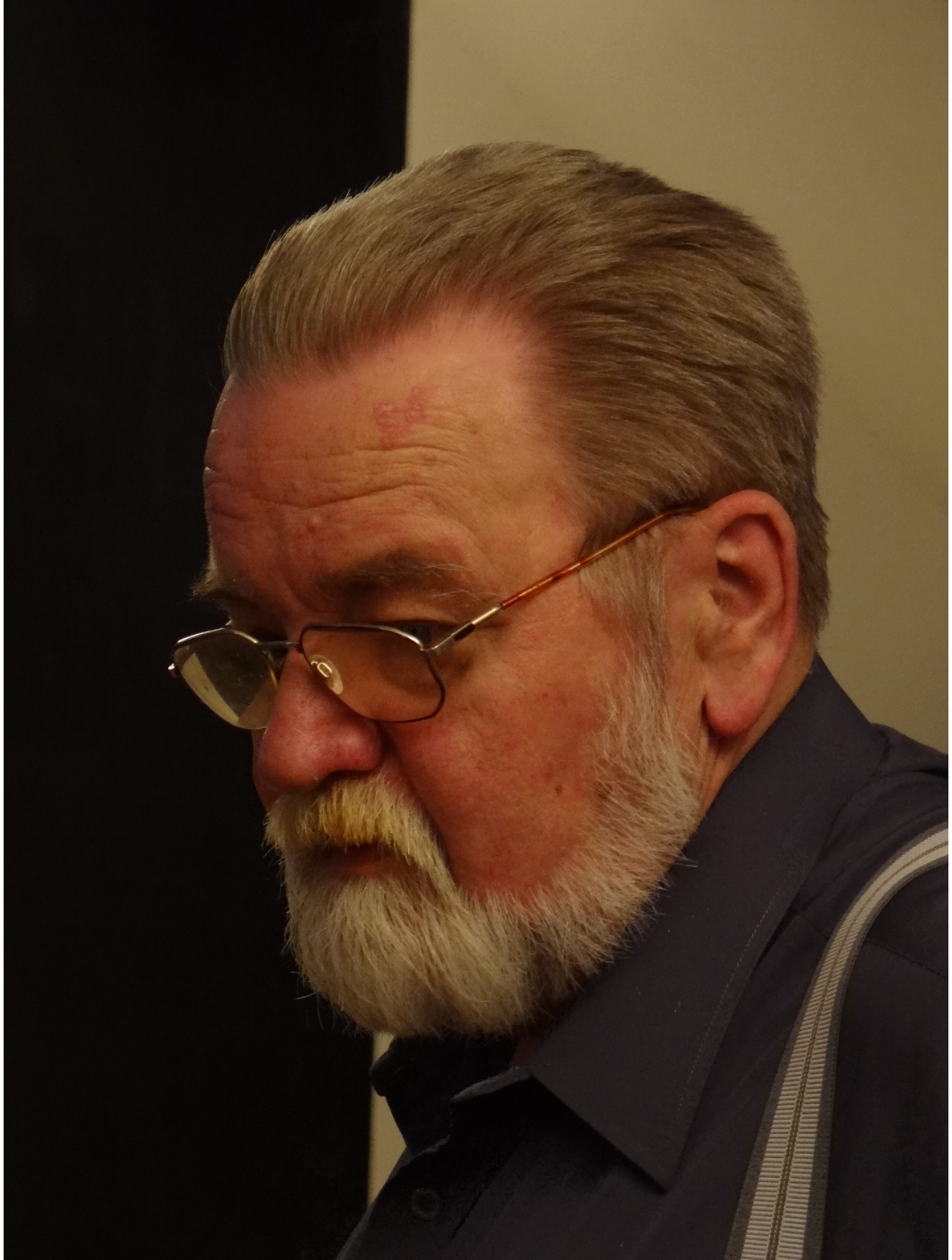}
\end{center}
\end{figure}

\section{Introduction}

The prominent Czech topologist Prof.~Petr Simon passed away 
on April~14th, 2018. 
He is an important link in the chain of renowned Czech topologists which 
includes well-known names like those of Eduard \v{C}ech,  
Miroslav Kat\v{e}tov, and Zden\v{e}k Frol\'{\i}k. 
In this note we wish to review his many achievements, where we concentrate
on his scientific contributions to the field of 
Set-Theoretic Topology.

Petr Simon was born on the 24th of February 1944 in Hradec Kr\'alov\'e in 
what is now the Czech Republic. 
He attended elementary school and secondary school in Prague, and studied at 
the Faculty of Mathematics and Physics of the Charles University in Prague 
in the period 1961--1966. 
Upon completing his studies he joined the Faculty as a Research Assistant 
in~1967; 
he worked there (mostly as a researcher) for the rest of his life. 
In~1977 he successfully defended his CSc 
(Candidate of Science, roughly equivalent to a PhD) 
dissertation 
\textsl{O lok\'aln\'{\i}m merotopick\'em characteru} 
(On local merotopic character) under the supervision of 
Professor Miroslav Kat\v{e}tov;
this came with an increase in rank to \emph{Researcher in mathematics}.
He was promoted to the rank of \emph{Independent researcher} in~1986 and,
after  defending his DrSc habilitation 
\textsl{Aplikace ultrafiltr{\accent23 u} v topologii} 
(Applications of ultrafilters in topology),  
to \emph{Leading researcher} in~1990. 
Finally, in~2001 he was awarded the title of full professor at 
Charles University.

\smallskip

When Petr Simon entered the mathematical life at the Charles University
there was a very active group of researchers in Topology, Set Theory and 
related areas. 
There was the \emph{Topological seminar} led by Miroslav Kat\v{e}tov, 
the \emph{Set Theory seminar}  of Petr Vop\v{e}nka,  
V\v{e}ra Trnkov\'a's \emph{Category Theory seminar},
and the seminars on \emph{Measure theory} and \emph{Uniform spaces} 
organized by Zden\v{e}k Frol\'{\i}k.
These were attended over the years by a talented group of
young mathematicians including 
Ji\v{r}\'{\i} Ad\'amek, 
Bohuslav Balcar, 
Lev Bukovsk\'y, 
Jan Hejcman,
Petr Holick\'y, 
Karel Hrb\'a\v{c}ek, 
Miroslav Hu\v{s}ek, 
V\'aclav Koubek,
Lud\v{e}k Ku\v{c}era,
V\v era K{\accent23 u}rkov\'a,
Vladim\'{\i}r M\"uller, 
Jaroslav Ne\v{s}et\v{r}il, 
Tom\'a\v{s} Jech, 
Jan Pachl, 
Jan Pelant, 
David Preiss, 
Karel P\v{r}\'{\i}kr\'y, 
Jan Reiterman, 
Vojt\v{e}ch R\" odl, 
Ji\v{r}\'{\i} Vil\'{\i}mov\-sk\'y,
Ji\v{r}\'{\i} Vin\'arek, and
Milo\v s Zahradn\'{\i}k, to mention but a few.

\smallskip

Every five years since 1961, Prague hosts the prestigious 
\textsl{Topological Symposium} an initiative of Eduard \v{C}ech. 
Petr Simon has participated in all but the first two of them, 
first as a speaker, and since~1981 as an organizer, and in~2006 as its 
chairman. 
He also participated actively in the organization  of the 
\emph{Winter School in Abstract Analysis, Section Topology}. 
These schools have been organized continuously since~1973 and were originally 
created by Zden\v{e}k Frol\'{\i}k as a winter getaway for the members of 
his seminars, but they quickly grew into 
important 
events in several fields (including, besides real and 
functional analysis, set theory and topology, also category theory and 
combinatorics). 
See~\cite{bs25WS} for a short history.

\smallskip
Petr Simon  was a longtime member of the editorial board of 
\emph{Topology and Its Applications} (from~1992 till~2018) and the 
managing editor of \emph{Acta Universitatis Carolinae Mathematica et Physica} 
(1989--2010). 
He was the representative of Czech Set Theory in the European Set Theory 
Society and the INFTY project (2009--2014).

\smallskip
To the best of our knowledge, Petr Simon was the advisor of the 
following doctoral students:
\begin{itemize}
\item Egbert Th\"ummel (1996):  
      \emph{Ramsey theorems and topological dynamics}
\item Eva Murtinov\'a (2002):  
       \emph{Separation Axioms in dense subsets}
\item Jana Fla\v{s}kov\'a (2006):  
      \emph{Ultrafilters and small sets}
\item David Chodounsk\'y (2011):  
      \emph{On the Katowice problem}
\item Jonathan Verner (2011):  
      \emph{Ultrafilters and Independent systems}
\item Jan Star\'y (2014):  
      \emph{Complete Boolean Algebras and Extremally Disconnected 
            Compact Spaces}
\end{itemize}

He has also supervised the Master's thesis of Dana Barto\v{s}ov\'a, and advised (albeit briefly) as PhD students Michael Hru\v{s}\'ak and Adam Barto\v{s}.

\smallskip
Petr Simon has published more than 70~research articles and participated in 
the preparation of the \emph{Handbook of Boolean Algebras}.
He co-edited the books  
\emph{The Mathematical Legacy of Eduard \v{C}ech} and  
\emph{Recent Progress in General Topology III}. 

His main contributions to topology lie in a thorough study of the 
combinatorial structure of the \emph{space of ultrafilters} $\betaN$
via almost disjoint and independent families, as well as the corresponding 
cardinal invariants of the continuum. 
His most lethal weapons were maximal almost disjoint (MAD) families 
and ultrafilters 
(\emph{the most beautiful things there are} --- he used to say).

\smallskip
He is and shall be missed by us and the whole community of Set Theory  
and Set-theoretic Topology. 
In what follows we take a more detailed look at some of his contributions 
to mathematics.

\textsc{Note to the reader}: we partitioned the references into two families.
The first, cited as~[S$m$], consists of works (co-)authored by Petr Simon;
the second, cited as~[O$n$], contains material related to Simon's work.

\section{Early beginnings}

Petr Simon published his first two articles in~1971 in the same volume 
of \emph{Commentationes Mathematicae Universitatis Carolinae} (CMUC), 
fittingly, 
the journal founded in~1960 by his academic grandfather 
Eduard \v{C}ech%
\footnote{The first volume of CMUC contains only two papers. 
One by Eduard \v{C}ech and the other by Miroslav Kat\v{e}tov.}.

In the first paper \cite{MR0288717} he made important contributions to the 
study of \emph{merotopic spaces}, 
introduced by Miroslav Kat\v{e}tov in~\cite{MR193608}, 
by John Isbell under the name \emph{quasi-uniform spaces} in~\cite{MR0170323} 
and again by Horst Herrlich in~\cite{MR0358706} as \emph{quasinearness spaces} 
--- a topic he revisited in \cite{MR0515013} and which formed part of his 
dissertation. 
He would later return to general continuity structures in his 
study of atoms in the lattice of uniformities and their relation to
ultrafilters \cites{MR0454928,MR0964058,MR0967749}, the latter two papers 
form joint work with Jan Pelant, Jan Reiterman and Vojt\v{e}ch R\"{o}dl, 
research undoubtedly stimulated by Zden\v{e}k Frol\'{\i}k's research 
seminar on Uniform Spaces.

In the second paper, \cite{MR0321008} (see also \cite{MR0355971}), 
he studied topological properties of Mary Ellen Rudin's Dowker space,
the then recently published first $\ZFC$ example of such 
a space~\cite{MR293583}.

\smallskip

In \cite{MR0339044} he contributed to the, then recently initiated, study 
of cardinal functions on topological spaces by proving that the cellularity 
of the square of a linearly ordered topological space  equals the density 
of the space itself.
This contains  Kurepa's result from \cite{MR0175792}
that for a linearly ordered space
the cellularity of the square is not larger than the successor of the 
cellularity of the space itself.

\smallskip

A joint paper with David Preiss \cite{MR0374875} contains a proof of the 
fact that a pseudocompact subspace of a Banach space equipped with the weak 
topology is compact. 
A key component of the proof extracts an important property of 
\emph{Eberlein compacta}%
\footnote{A topological space is Eberlein compact if it is homeomorphic to 
          a weakly compact subset of a Banach space.} 
later dubbed the Preiss-Simon property: 
a space $X$ has the \emph{Preiss-Simon property} if for every 
closed $F\subseteq X$, each point $x\in F$ is a limit of a sequence of 
non-empty open subsets of~$F$. 
He continued the study of Eberlein compact spaces in~\cite{MR0407575}, 
providing partial results toward proving that continuous images of 
Eberlein compacta are Eberlein;
this was proved soon thereafter by 
Yoav Benyamini, Mary Ellen Rudin and Michael Wage in~\cite{MR625889}. 

\smallskip

One of Petr Simon's strengths was his ability to construct ingenious 
examples and counterexamples in topology. 
Answering a question raised by J. Pelham Thomas in~\cite{MR0236884}
he constructed  in~\cite{MR0487978} the first example of an infinite 
maximal connected Hausdorff space, that is,
a connected Hausdorff topological 
space~$X$ such that any finer topology on~$X$ is disconnected. 
In \cites{MR0517197,MR0588856} he constructed a compactification~$b\N$ 
of the countable discrete space~$\N$ with a sequentially compact 
remainder such that no sequence in~$\N$ converges in~$b\N$. 
For later examples, one can look at his construction in~\cite{MR1893957} of 
a connected metric space in which every infinite separable subspace is not 
connected, or a joint work with Steve Watson \cite{MR1183326} where a 
completely regular space which is connected, locally connected and 
countable dense homogeneous, but not strongly locally homogeneous is 
constructed.

\section{Almost disjoint families and ultrafilters}
\label{sec:adf+uf}

The \emph{Nov\'ak number} $\enn(X)$%
\footnote{The \emph{Baire number} for the rest of the world, 
          and eventually even for the Czech community.} 
of a topological space without isolated points is the minimum cardinality of 
a family of nowhere dense subsets that covers~it.  
Petr Simon first used the Nov\'ak number in \cites{MR0550459,MR0580656} 
to give alternative proofs of a result of Mikhail Tkachenko.

The study of Nov\'ak number was a popular research topic in Prague. 
For example, Petr \v{S}t\v{e}p\'anek and Petr Vop\v{e}nka showed 
in~\cite{MR225657} that the Nov\'ak number of any nowhere separable 
metric space is at most~$\aleph_1$. 
Petr Simon \cite{MR0464133} greatly generalized their result, and  then 
together with Bohuslav Balcar and Jan Pelant he studied the behaviour of 
the Nov\'ak number of the \v{C}ech-Stone remainder~$\Nstar=\betaNminN$ 
of~$\N$  in~\cite{MR0600576}. 
In this important  paper, they showed that the value of $\enn(\Nstar)$ depends 
heavily on set-theoretic axioms.  
To study these they introduced the \emph{distributivity number} $\haa$ 
of the Boolean algebra $\pow(\N)/\fin$ and proved their celebrated and 
influential \emph{Base Tree Theorem} which states that the algebra has a 
dense subset~$T$ which is a~$\cee$-branching tree of height $\haa$, 
and showed that 
(1) if $\haa<\cee$ then $\haa\leq \enn(\Nstar)\leq \haa^+$, 
(2) if $\haa=\cee$ then $\haa\leq \enn(\Nstar)\leq 2^\cee$, and 
(3) $\haa=\enn(\Nstar)$ if and only if there is a Base Tree without 
    cofinal branches.

\smallskip

Simon quickly put the theorem to further use in \cite{MR0493972} to show 
that if $\enn(\Nstar )> \cee$ then $\Nstar$~contains a dense linearly ordered 
subspace, in~\cite{MR1346339} where it is shown that $\haa$~is the minimal 
size of a family of sequentially compact spaces whose product is not 
sequentially compact and in~\cite{MR1419359} where he constructs 
a $\sigma$-centered atomic, almost rigid tree-like Boolean algebra such 
that every injective endomorphism is onto, and every surjective homomorphism 
is injective.

\smallskip

The paper \cite{MR0464133} was also the start of a long and fruitful 
collaboration of Petr Simon with Bohuslav Balcar. 
The two have co-authored over twenty publications and co-directed the 
well-known 
\emph{Semin\'a\v{r} z po\v{c}t\accent23 u}%
\footnote{The name is a play on the curious fact that, at the time, set 
          theory was being introduced in mathematics classes taught in lower 
          elementary school, these classes were then called ``reckoning''} 
(loosely translated as ``Seminar on reckoning'' or ``Seminar on counting'') 
and together raised a new generation of Czech Set-Theorists and Topologists.

\smallskip

Together with Peter Vojt\'a\v{s} \cites{MR0621987,MR0628639} 
(see also \cites{MR0818228,MR0991597}) they continued 
the project initiated by Balcar and Vojt\'a\v{s} in determining which 
subsets of Boolean algebras admit a disjoint refinement. 
To discuss a particularly important case of the above mentioned phenomenon 
let us introduce some notation first: 
A family~$\calA$ of infinite subsets of~$\N$ is 
\emph{almost disjoint} (AD) if every two distinct elements of~$\calA$ 
have  finite intersection. 
An infinite AD family is \emph{maximal} (MAD) if it is maximal with respect to 
inclusion or, equivalently, if for every infinite $B\subseteq\N$ there is 
an~$A\in\calA$  such that $B\cap A$ is infinite. 
Given a MAD family~$\calA$ we denote by~$\calI(\calA)$ 
the family of $\calA$-small subsets of~$\N$: 
those sets which have infinite intersection with only finitely many 
elements of~$\calA$, and by~$\calI^+(\calA)$ the family 
of those subsets of $\N$ which are not $\calA$-small, and therefore
are called \emph{$\calA$-large}.

The crucial problem, first studied in \cite{MR0818228}, 
is~$\mathrm{RPC}(\omega)$: 
\begin{quote}
Given a maximal almost disjoint family  $\calA$, does there exist an 
almost disjoint refinement for the family of all $\calA$-large sets?
\end{quote}
A closely related problem was formulated independently and in a different 
context by Paul Erd\H{o}s and Saharon Shelah in~\cite{MR319770}: 
\begin{quote}
Is there a completely separable MAD family, that is,
is there a MAD family $\calA$ that is itself an almost disjoint refinement 
for the family of $\calA$-large sets?
\end{quote}

The problem $\mathrm{RPC}(\omega)$ has the following two topological 
equivalents
\begin{enumerate}
\item for every nowhere dense subset $N$ of~$\Nstar$ there is 
      a family of $\cee$-many disjoint open subsets of~$\Nstar$ 
      each having~$N$ in its closure, and 
\item for every nowhere dense subset~$N$ of~$\Nstar$ there 
      is a nowhere dense subset~$M$ of~$\Nstar$  such that $N$~is a 
      nowhere dense subset of~$M$, see~\cite{MR1056181}. 
\end{enumerate}
It is believed that both problems have positive solutions in~$\ZFC$, and Simon 
and collaborators have made steady progress toward a solution 
\cites{MR0818228, MR0991597, MR1600453}. 
Currently the best result is due to Saharon Shelah who showed 
in~\cite{MR2894445} that the answer is positive if $\cee<\aleph_\omega$ 
and that a negative solution would imply consistency of the existence of 
large cardinals.

\smallskip
Having completed their study of the structure of the Boolean algebra 
$\pow(\omega)/\fin$, 
Simon and Balcar turned to the higher cardinal analogues  
$\Pkappamodsmall$,  
a study initiated by Balcar and Vop\v{e}nka in~\cite{MR314618} 
and the corresponding study of \emph{uniform ultrafilters on $\kappa$}. 
In \cites{MR0982781,MR0991597}, extending earlier results of Balcar 
and Vop\v{e}nka, they showed that $\haa_\kappa=\aleph_0$ for cardinals $\kappa$ 
of uncountable cofinality and $\haa_\kappa=\aleph_1$ for uncountable 
cardinals $\kappa$ of countable cofinality%
\footnote{Here we take as definition of $\haa_\kappa$ the minimal~$\mu$ 
          such that forcing with $\Pkappamodsmall$~adds a new 
          set of size~$\mu$.}, 
while in~\cite{MR1357756} they computed the Nov\'ak numbers of the 
space $U(\kappa)$ of \emph{uniform ultrafilters on $\kappa$}, 
i.e., the Stone space of the Boolean algebra 
$\Pkappamodsmall$. 
They proved that 
(1) $\enn(U(\kappa))=\aleph_1$ for all~$\kappa$ of uncountable cofinality, and 
(2) $\enn(U(\kappa))=\aleph_2$ for all~$\kappa$ of countable cofinality, 
assuming one of~$\neg\CH$, 
$2^{\aleph_1}=\aleph_2$, and
$\kappa^{\aleph_0}=2^\kappa$ holds. 
They further studied which cardinals are being collapsed in generic extensions 
by  $\Pkappamodsmall$ in 
\cites{MR0982781,MR0991597,MR1835243}; 
the final word on this is yet to be written. 
Currently the best partial results, by Saharon Shelah, appear 
in~\cite{MR2298480}.

\smallskip

Continuing the work of Balcar and Vojt\'a\v{s}, 
they further studied the structure of $U(\kappa)$ along the lines 
of~$\mathrm{RPC}$, focusing on a question of 
Wis Comfort and Neil Hindman from~\cite{MR429573}: 
\emph{Given an arbitrary cardinal $\kappa$, is each point 
in~$U(\kappa)$ a $\kappa^+$-point}, i.e., is there a family of $\kappa^+$  
pairwise disjoint open subsets of~$U(\kappa)$ each containing the point 
in its closure? 
They showed in~\cite{MR0673783} that the answer is positive for regular 
cardinals and later Simon showed in~\cite{MR0778069} that it
is also positive for cardinals of countable cofinality. 
For such cardinals~$\kappa$, answering another question of Wis Comfort, he 
showed in \cite{MR2103088} that there exist uniform ultrafilters on~$\kappa$
that cannot be obtained from the set of all sub-uniform ultrafilters 
by iterating the closure of sets of size less than~$\kappa$.

\smallskip

In a pair of articles \cites{MR1129703,MR1189823} Balcar and Simon 
investigated the minimal $\pi$-character of ultrafilters on Boolean algebras. 
They showed that if a Boolean algebra is homogeneous or complete then the 
minimal $\pi$-character coincides with the \emph{reaping number} of the 
Boolean algebra, which is defined as the minimal size of a family of 
non-zero elements such that no element splits them all. 
Then they used their results to show that every extremely disconnected 
ccc space in which $\pi$-weight and minimal $\pi$-character coincide, 
contains a point which is not an accumulation point of a countable discrete 
set, also known as a \emph{discretely untouchable point}
(see~\cite{MR1108207}).

To put this result  in context, we recall that Zden\v{e}k Frol\'{\i}k 
proved in~\cite{MR264584} that every infinite compact extremely 
disconnected space (the Stone space of a complete Boolean algebra) 
is not homogeneous.
This proof, however, did not produce examples of simple topological properties
shared by some but not all of the points of the space. 
It is clear that every compact space contains points which are accumulation 
points of countable discrete sets, so ``being discretely untouchable''
is a property not shared by all points in a compact space ---
if it were shown to be shared by some then one would have a more concrete
reason for non-homogeneity. 

The non-homogeneity of~$\betaNminN$ was proved under~$\CH$ by Walter Rudin 
in~\cite{MR80902} 
using the simple property of being a \emph{P-point} 
(the intersection of countably many neighbourhoods is again a neighbourhood),
in~$\ZFC$ by Zden\v{e}k Frol\'ik in~\cite{MR266160} using a rather 
unwieldy invariant of points in~$\betaNminN$ that involved copies of the 
space inside itself,
and again in~$\ZFC$ by Ken Kunen in~\cite{MR588822} using the property
of being a \emph{weak P-point} 
(not an accumulation point of any countable set).
This latter proof introduced a powerful method of constructing ultrafilters 
using independent linked families.%
\footnote{The example of such a family in Kunen's paper is given by an explicit
          formula due to Petr Simon, rather than Kunen's own construction
          ``involving trees of trees''} 
This method was used by Simon in~\cite{MR0824025} and~\cite{MR0863940} 
to construct ultrafilters on~$\N$ that are minimal in the Rudin-Frol\'{\i}k 
order, and in~\cites{Simon1984,MR0869226} to give a $\ZFC$~construction 
of a closed separable subspace of~$\Nstar$ that is not a retract of~$\betaN$.

That last construction is highly involved and takes about ten~pages.
It also illustrates what some us know well, namely that Petr Simon had a way 
with words; after an outline of the construction and just before the hard 
work starts we find the following gem:
``so let us begin now to swallow the indigestible technicalities''.

\smallskip

In joint work Petr Simon, Murray Bell and Leonid Shapiro \cite{MR1328339} 
showed that a large class of spaces,
called \emph{orthogonal $\Nstar$-images},
that includes 
all separable compact spaces, 
all compact spaces of weight at most~$\aleph_1$ and 
all perfectly normal compact spaces, 
is closed under products of $\cee$-sized families.
A compact space~$K$ is an orthogonal $\Nstar$-image if $\Nstar\times K$
is a continuous image of~$\Nstar$. 
This result --- trivially true under~$\CH$ by Parovichenko's Theorem ---  
puts some limits to the \emph{rigidity phenomena} that one encounters 
assuming~$\PFA$, see, e.g., Farah~\cite{MR2143881} 
and Dow-Hart~\cite{MR1679586}.

\smallskip

The Stone duality that links zero-dimensional compact Hausdorff spaces 
and Boolean algebras permeates much of Petr Simon's work. 
In a joint paper with Martin Weese \cite{MR0803920} he constructed
non-homeomorphic \emph{thin-tall} scattered spaces, later superseded by 
his work with Alan Dow in~\cite{MR1157434}.  
In~\cite{MR1921660}, answering a question of Alexander Arkhangel'ski\u{\i}, 
he showed that there is a separable compact Hausdorff space~$X$ having 
a countable dense-in-itself, dense set~$D$ which is a P-set in~$X$.

\smallskip

Another constant in Simon's work is the use of cardinal invariants of 
the continuum, both as tools \cites{MR0961403,MR1056363} and as principal 
objects of study~\cites{MR1263799,MR2138279,MR2361889,MR3072778}.
For example, the first of these four paper shows that Sacks forcing collapses
$\cee$ to~$\bee$.

\smallskip
In a joint paper with Alan Dow and Jerry Vaughan \cite{MR0961403} he found 
one of the first applications of Set Theory in Algebraic Topology and 
Homological Algebra. 
The study of the interactions between these fields has started to flourish 
only recently.


\section{Convergence properties}

Petr Simon made a strong impact on the study of Fr\'echet spaces and 
their generalizations, often in collaboration with his Italian colleagues 
Angelo Bella, Camilo Costantini and Gino Tironi. 
Recall that a topological space $X$ is \emph{Fr\'echet} if for every point 
in the closure of a set $A\subseteq X$ there is a sequence of elements of~$A$ 
that converges toconverging to~it. 
There is an extremely close connection between Fr\'echet spaces and almost 
disjoint families.

\smallskip

In \cite{MR0597764} Simon gave the first $\ZFC$ example of a compact Fr\'echet
space whose square is not Fr\'echet. 
The result followed using known techniques from the fact that there is a MAD 
family~$\calA$ that can be partitioned into two subfamilies, 
each of which is not maximal when restricted to an $\calA$-large set,
they are thus called \emph{nowhere maximal}. 
Simon's proof of this is a thing of beauty, a \emph{proof from the book}, 
a short concise proof by contradiction. 
He later also used this result in~\cite{MR0749120}.

In \cite{MR1696588} and \cite{MR1922133} Simon studied the question of 
Tsugunori Nogura, whether the product of two Fr\'echet spaces neither 
of which contains containing a copy of the sequential fan can contain 
a copy of it. 
First, \cite{MR1696588} he showed that assuming $\CH$ such spaces can be 
constructed, and later in joint work with Tironi, \cite{MR1922133}, gave 
partial results in the opposite direction. 
The final solution to the problem was given by Stevo Todor\v{c}evi\'c 
in~\cite{MR1999941}
extending the approach of Simon and Tironi:
the Open Colouring Axiom implies a positive answer.
 
On the other hand, Costantini and Simon \cite{MR1783423} answering 
another question of Nogura gave a $\ZFC$ example of two Fr\'echet spaces, 
the product of which does not contain a copy of the sequential fan but fails 
to be Fr\'echet. 
For this they used a construction of an AD family with properties resembling 
those of a completely separable MAD family. 
Later Bella, Costantini and Simon \cite{MR2227019} showed that assuming~$\CH$ 
one can construct Fr\'echet spaces containing a copy of the sequential fan 
which are pseudocompact. 
In his last paper concerning the subject \cite{MR2419371} Simon further 
advanced these techniques to construct a Fr\'echet space not containing 
a copy of the sequential fan all of whose finite powers are Fr\'echet.

\smallskip

By allowing transfinite sequences when reaching points in the closure one 
arrives at a more general notion of a \emph{radial space} or 
a \emph{Fr\'echet-chain-net space}. 
If one only requires that every non-closed set contains a transfinite 
sequence converging outside of the set, one gets the notion of 
a  \emph{pseudoradial space} or a \emph{chain-net space}, just as with 
the usual convergence one defines a \emph{sequential space}. 
Simon started to study such spaces in a joint paper~\cite{MR0647028} with 
Ignacio Jan\'e, Paul Meyer and Richard Wilson  where they, 
assuming~$\CH$, constructed Hausdorff examples of countably tight pseudoradial 
spaces which are not sequential. 
In a paper with Tironi \cite{MR0843427} they produced a $\ZFC$ example. 
Completely regular examples were given soon after by Juh\'asz and Weiss 
in~\cite{MR874661}.
In \cites{MR1378621,MR1679197} and \cite{MR1811550} Simon and co-authors 
look at products of (compact) pseudo-radial spaces.

\smallskip

A further weakening gives rise to Whyburn and  weakly Whyburn spaces: 
We say $X$~is a \emph{Whyburn} space if whenever $x\in \overline A\setminus A$, 
there is a $B\subseteq A$ such that $ \overline B\setminus A=\{x\}$, 
and $X$~is \emph{weakly Whyburn} if whenever $A\subseteq X$ is not closed, 
there is a $B\subseteq A$ such that $|\overline B\setminus A|=1$. 
Simon first studied these spaces in~\cite{MR1307264} 
(using different terminology) and showed that there are two Whyburn 
spaces whose product is not weakly Whyburn.
In \cite{MR2352741}, assuming~$\CH$, Bella and Simon construct a 
pseudocompact Whyburn space of countable tightness that is not Fr\'echet.

\smallskip

The paper \cite{MR2024948} of Bella and Simon continues the study 
initiated in~\cite{MR1918105} by Dow, Tkachenko, Tkachuk and Wilson
of \emph{discretely generated} spaces: spaces where points in the 
closure can be reached by discrete sets --- another weakening of radiality. 
They show that countably tight countably compact  spaces are discretely 
generated and show that this consistently fails for pseudo-compact spaces 
of countable tightness. 

\smallskip

The paper \cite{MR2024948} also contains results concerning spaces of 
continuous functions endowed with the topology of pointwise convergence. 
In particular, it shows that if $X$ is $\sigma$-compact then $C_p(X)$~is 
discretely generated.
The work with Tsaban \cite{MR2361889} shows that the 
\emph{pseudo-intersection number $\pee$} is the minimal cardinality 
of a set~$X$ of reals, such that $C_p(X)$~does not have the Pytkeev property.
This is another local property of topological spaces: 
A space~$X$ has the \emph{Pytkeev property} if for every $A\subseteq X$  
and every $y\in\overline A\setminus A$ there is a countable
family~$\calA$ of infinite subsets of~$A$ such that every 
neighbourhood of~$x$ contains a member of~$\calA$. 
 
The spaces of continuous functions over the Mr\'owka-Isbell spaces associated 
to almost disjoint families in relation to the study of the Lindel\"of 
property in~$C_p(X)$ are investigated in~\cite{MR2238729}.

\smallskip

Function spaces with the topology of uniform convergence were considered 
by Bella and Simon in \cite{MR1047521} where it is shown that the set of 
nowhere constant functions is dense in~$C(X,Y)$  if $Y$~is a normed linear 
space and $X$ is a dense-in-itself normal space  or separable completely 
regular space.

\section{Other}

In his undergraduate topology classes Petr Simon would maintain that 
\emph{every respectable topological space is Tychonoff}. 
Even so, he himself has sinned and occasionally looked at the less 
respectable ones. 
We already saw an example of this above 
\cites{MR0647028, MR0843427}. 
In \cite{MR0749120}, using his \emph{partitionable MAD family} and 
the so-called \emph{Jones machine} he produced an example of two regular, 
functionally Hausdorff spaces such that the product of their completely 
regular modifications does not coincide with the completely regular 
modification of their product%
\footnote{The \emph{completely regular modification} of a space $(X,\tau)$ 
          is $X$ endowed with the weakest topology making all 
          $\tau$-continuous functions continuous}. 
In a joint paper \cite{MR1082937} with Eraldo Giuli he showed that the 
category of all topological spaces in which every bounded set is Hausdorff 
is not co-well-powered.
In \cite{MR2054805}, he together with Gino Tironi showed that locally 
feebly compact first countable regular spaces can be densely embedded 
into feebly compact first countable regular spaces. 
In Bela, Costantini and Simon \cite{MR2227019} appears a consistent 
construction of a countably compact Hausdorff space which is Fr\'echet 
and contains a copy of the sequential fan.

\smallskip

In a joint paper with Pelant and Vaughan \cite{MR0978716} we find information
about the minimal number of free prime filters of closed sets on a 
non-compact space.
For completely regular spaces this number is at least~$\aleph_2$, 
while for Hausdorff spaces the best bound they could find is~$\aleph_1$.  

Simon together with Hindman and van Mill \cite{MR1089407} considered 
$\beta\Z$ as a compact left topological semigroup, and show that 
there is a strictly increasing chain of principal left ideals and of 
principal closed ideals.

\smallskip
The  paper \cite{MR2282715} with Fred Galvin, which answers a 1947 problem 
of Eduard \v{C}ech by constructing  a so called \emph{\v{C}ech  function} --- 
a pathological closure operator on $\pow(\omega)$ which is surjective 
yet not the identity, has a curious history: Galvin knew since 1987 that
the existence of a completely separable MAD family suffices, 
while Balcar, Do\v{c}k\'alkov\'a and Simon in \cite{MR0818228} (1984) 
constructed an AD family with similar properties, which would already suffice 
for the Galvin result. 
\begin{figure}[ht]
\begin{center}
 \includegraphics[width=.7\hsize]{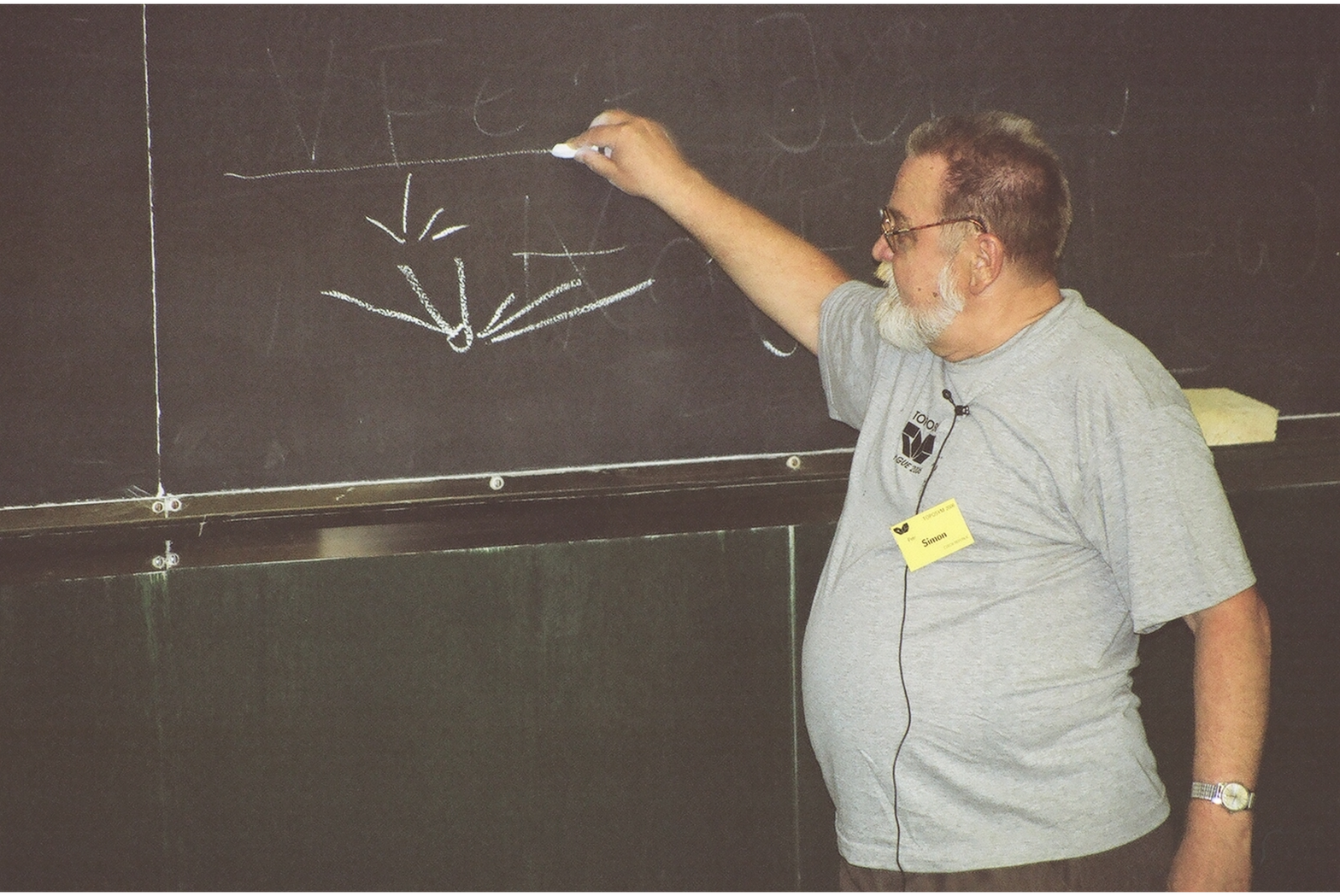}
\end{center}
{\small Petr Simon lecturing on the \v{C}ech function during TOPOSYM 2006.}
\end{figure}
Even though they were both well aware of each others results, it took 
them 20 more years to put the two facts together.

\smallskip

Convergence structures on groups were treated in \cite{MR0904772} 
and \cite{MR1371690}. 
In the first paper Simon and Fabio Zanolin show that there is a 
Boolean coarse convergence group that cannot be  embedded into a sequentially 
compact convergence group, while in second one Simon uses the 
additive group of rational numbers to show that the theory of sequential 
groups does not admit a reasonable notion of completeness. 
To exemplify this he shows that the smallest convergence structure  
on~$\Q$ making the sequence~$\langle \frac 1n\rangle_n$ converge to~$0$ 
is complete, and constructs  another group convergence on~$\Q$ such that some 
irrationals, but not all, are limits of Cauchy sequences.

\smallskip

Apart from his research activities, 
Petr Simon has written several surveys 
and introductory articles 
\cites{MR0624782, MR0645299, MR0991618, MR1269145, MR1617557, MR2380254} 
as well as bibliographical articles dedicated to the life and mathematics 
of Bohuslav Balcar~\cite{MR3914709}, 
Eduard \v{C}ech \cites{MR1241376, MR1249622, MR1280444,MR1282960} 
(including co-editing the book \emph{The mathematical legacy of Eduard \v{C}ech}
  with Miroslav Kat\v{e}tov \cite{cech-book}), 
Zden\v{e}k Frol\'{\i}k \cites{ MR1071064,MR1084905,MR1222259}, 
Miroslav Kat\v{e}tov \cites{MR1408306,MR1446403,MR1617561} and 
Jan Pelant \cites{MR2206282,MR2211007}. This all, of course, reflected his 
status in the mathematical community in Prague.

\section{Questions}

The papers co-authored by Petr Simon contain many questions.
From our conversations with him we got the impression he would very much
have liked to be instrumental in solving the following four.

\begin{question}[\cite{MR0869226}]
Must a closed separable subset of $\betaN$ that is not a retract
have a tiny sequence?
\end{question}

A \emph{tiny sequence} in a space~$X$ is a family $\{U_{n,m}: n,m\in \omega\}$ 
(a matrix) of open sets with the property that
for the union $\bigcup_{m\in\omega}  U_{n,m}$ of each row is dense,
and yet whenever one chooses finite subfamilies $\{U_{n,m}:m\in F_n\}$
of the rows the union $\bigcup_{n\in\omega}  \bigcup_{m\in F_n}U_{n,m}$ is not dense.
These were introduced by Szyma\'nski in~\cite{MR744396} and used there in 
the proof that a particular separable subset of~$\betaN$,
constructed from Martin's Axiom, was not a retract. 
Simon's $\ZFC$-example has a natural tiny sequence, hence the question. 

In~\cite{MR810825} Leonid Shapiro also constructed a separable non-retract
of~$\betaN$.
The construction is indirect: after deriving a condition under which the
absolute of a compact space is not an absolute retract of~$\betaN$ he
constructs a separable compact space of weight~$\aleph_1$ that does not meet
that condition.
We do not know whether this yields in fact a counterexample to this question.

\begin{question}[\cite{MR0991597}]
$\mathrm{RPC}(\omega)$: 
Given a maximal almost disjoint family~$\calA$ does there exist 
an almost disjoint refinement for the family of all $\calA$-large sets?
\end{question}

We discussed this question in Section~\ref{sec:adf+uf}
and all we can say is that it is a beautiful question, both in its
combinatorial formulation and its topological equivalents.

\begin{question}[\cite{MR1056363}]
Is there (in $\ZFC$) a non-meager P-filter? 
\end{question}

We mentioned that Walter Rudin's non-homogeneity proof for~$\betaNminN$ 
from~\cite{MR80902} used P-points. 
These points feature, as ultrafilters, in solutions to many combinatorial 
problems involving the set of natural numbers and quite often they are even 
necessary for that solution.
When Saharon Shelah showed that the existence of P-points in~$\betaNminN$
is not provable in~$\ZFC$ alone, see~\cite{MR728877}, this spurred research 
into good approximations of P-points whose existence could be shown on 
basis of~$\ZFC$ alone. 
In the paper \cite{MR1056363} Petr Simon and his co-authors took up this topic. 
They focused on the question of how close, in~$\ZFC$, a P-filter can be to 
an ultrafilter and hypothesized that it can be non-meager%
\footnote{As a subset of the Cantor space; since filters here are 
          collections of subsets of the natural numbers, we can think of 
          them as subsets of the Cantor space.}. 
That paper contains Petr Simon's proof of the existence of such filters 
under $\tee=\bee$ or $\bee<\dee$ and shows that their non-existence would 
need large cardinals. 
Whether they can be shown to exist in $\ZFC$ is still an intriguing 
open question.

\begin{question}[\cite{MR1108207}]
Does every extremally disconnected compact space contain a 
        discretely untouchable point?
\end{question}

This was also discussed in Section~\ref{sec:adf+uf}. 
As mentioned there this is an attempt to prove that compact extremally
disconnected are not homogeneous by means of an easily formulated property 
`shared by some but not all points'.
Other properties have been brought to bear on this problem but 
`discretely untouchable' seems particularly susceptible to combinatorial
treatment.

\section*{Acknowledgments}

We thank Kl\'ara Chytr\'a\v{c}kov\'a for supplying us with many biographical
details about her father Petr Simon. We also wish to thank  Miroslav Hu\v sek 
and Jaroslav Ne\v set\v ril for commenting on an early draft of the paper, 
thus helping us to make it more historically accurate.

\end{document}